\documentclass[12pt,reqno]{article}

\title {Uniform Infinite Planar Triangulations}
\author{Omer Angel \and Oded Schramm}

\date{February, 2003}  % Revision

\usepackage{latexsym, amsfonts, amsmath, amsthm, graphicx}

\numberwithin{equation}{section}
\numberwithin{figure}{section}

\newtheorem{thm}{Theorem}[section]
\newtheorem{prop}[thm]{Proposition}
\newtheorem{coro}[thm]{Corollary}
\newtheorem{conj}[thm]{Conjecture}
\newtheorem{lemma}[thm]{Lemma}

\newtheorem{defn}[thm]{Definition}

\newcommand{\N}{{\mathbb N}}
\newcommand{\Z}{{\mathbb Z}}
\newcommand{\R}{{\mathbb R}}
\renewcommand{\P}{{\mathbb P}}
\newcommand{\T}{{\mathcal T}}

\renewcommand{\phi}{\varphi}
\newcommand{\ep}{\varepsilon}

\newcommand{\note}{{\noindent \bf Note.\ }}

\newcommand{\be}{\begin{enumerate}}
\newcommand{\ee}{\end{enumerate}}
\newcommand{\bi}{\begin{itemize}}
\newcommand{\ei}{\end{itemize}}

\begin{document}
\maketitle

\begin{abstract}
The existence of the weak limit as $n\to\infty$ of the uniform measure
on rooted triangulations of the sphere with $n$ vertices is proved. Some
properties of the limit are studied. In particular, the limit is a
probability measure on random triangulations of the plane.
\end{abstract}


\section{Introduction}    %\input intro.tex

%%%%%%%%%%%%%%%%%%%%%%%
\subsection{Motivation}
%%%%%%%%%%%%%%%%%%%%%%%

What is a generic planar geometry?

There are many different planar geometries. The most commonly used one is the
Euclidean plane, but is it generic? Is it more natural than, say, the
hyperbolic plane?

For simplicity, consider discrete planar geometries (realized as planar
graphs). Now there are still many choices. The lattice $\Z^2$ is the graph
most commonly associated with planar geometry, but there is no a priori
reason to prefer it over the triangular lattice, or any other lattice. One
possible approach is based on convenience, preferring at each time the most
convenient framework to work with. Even by that criterion no single geometry
is always the best. Thus, some recent results are naturally adapted to
the triangular lattice \cite{Smir}.

When we use a lattice, we force much more structure into our geometry than
the topological condition of planarity necessitates. Random planar graphs,
such as Delaunay triangulations, have less enforced structure, but they still 
arise from the underlying Euclidean geometry. Is there a clear reason to
prefer the Euclidean over the hyperbolic plane?

The approach used here is to consider a probability measure --- in some
sense a uniform measure --- on planar geometries. Then we can ask what
properties does a typical sample of that measure have. The way this is done
is by considering discrete geometries, realized in the form of infinite
planar triangulations, and finding an interesting distribution on them.
Over finite planar triangulations the uniform measure is a natural choice.
We prove the existence of a probability measure on infinite planar
triangulations which is the limit of the uniform distributions on finite
planar triangulations as their size tends to infinity. A sample of this
measure is called the uniform infinite planar triangulation (UIPT). This
model was suggested in \cite{BeSc}, where Benjamini and Schramm show a.s.\ 
parabolicity of a wide class of distributions on infinite planar graphs
under the condition of a uniform bound on the vertex degrees. Alas, the
results there require the vertex degrees to be bounded, and
hence do not apply to the UIPT.

\medskip

The uniform {\em finite} planar triangulation and related objects have been
studied by both combinatorists and physicists. Mathematical study is
traced back to the 1960's with Tutte's attempts at the four color problem.
In a series of papers Tutte was able to count the number of planar maps of
a given size of various classes, including triangulations
\cite{census,census1,census2,census3}. One of the conjectures he raised is 
that almost all planar maps are asymmetric, i.e., have no non-trivial
automorphisms.

Tutte later proved his conjecture for a specific class of planar maps
\cite{Tut2}. Random planar maps (and triangulations among them) have been
studied extensively since then by others, proving Tutte's conjecture in a
more general setting \cite{RiWo2}.

Previous research here focused on finite triangulations, but many of the
results are about the asymptotic properties of planar maps and can be
translated directly into claims about the infinite triangulations we study.
Thus, there are results about the distribution of degrees in a uniformly
chosen triangulation \cite{GaRi}, the size of 3-connected components
\cite{BRW}, and probabilistic 0-1 laws \cite{BCR}. A key novel feature of
this paper is that we consider instead a distribution on {\em infinite}
maps. As it turns out, in some respects it is easier to work with the
infinite object than it is with the finite maps.

Schaeffer found a bijection between certain types planar maps and labeled
trees \cite{Scha}. Chassaing and Schaeffer~\cite{ChSc} recently used that
bijection to show a connection between the asymptotic distribution of the
radius of a random map and the integrated super-Brownian excursion. They
deduce from this connection that the diameter of such a map of size $n$ 
scales as $n^{1/4}$. While they work with planar quadrangulations and we 
with triangulations, it appears that such local differences are
insignificant when large scale observations such as diameter, growth,
separation, etc. are concerned. This phenomenon is referred to as
universality.

The physicists study such triangulations under the titles of {\em dynamic
triangulations}, closely related to {\em 2-dimensional quantum gravity}.
There, the essential idea is to develop a quantum
theory of gravity by extending to higher dimensions the concept of Feynman
integrals on paths. Triangulations are used as a discretized version of a 2 
dimensional manifold, and a function is averaged over all of them
\cite{AmWa, BoKa}. Physicists are more interested in a continuous scaling
limit of the discrete model, which is believed to exist.

Physicists introduced here the methods of random matrix models \cite{FGZJ}.
Through these methods and other heuristics many conjectures were made on
the structure of such triangulations. In particular, it is believed 
that the Hausdorff dimension of the scaling limit of 2-dimensional quantum
gravity is 4 \cite{AmWa}. For a good general exposition of quantum gravity
see \cite{ADJ}, as well as \cite{Amb,Dav}.

Of particular interest
is the KPZ relation \cite{KPZ} which relates
critical exponents for a number of models on the plane and in 2 dimensional 
quantum gravity. This relation has been used to predict various exponents
such as non intersection exponents for Brownian motion in the plane
\cite{KPZ1,KPZ2}. Later a rigorous derivation of the same values was found
using the $SLE$ process \cite{LSW1,LSW2,LSW3}.

It is hoped that this work will be the foundation for a rigorous study into
the scaling limit of random planar maps, and thereby enable a better
understanding of the relation between random surfaces and critical models
on smooth surfaces.

\medskip

Section~\ref{sec:count} summarizes some results on counting triangulations
which are the basis for much of what follows. Section~\ref{sec:basic}
describes some properties of the UIPT that follow directly from the
formulas for counting triangulations. In particular, is shown that a.s. the 
UIPT has one end, i.e., the limiting process does not add any topological
complications to the triangulation.

In Section~\ref{sec:limit} the existence of the limit distribution is
proved. In Section~\ref{sec:local} we gives another characterization of the
UIPT by a locality property. This roughly means that different regions in
the triangulation are independent of one another and that each region is
uniformly distributed among all triangulations of a given size (and hence
the name uniform triangulation for the infinite graph).
Section~\ref{sec:GW} describes a multi-type Galton-Watson tree naturally
associated with a UIPT.

In Section~\ref{sec:equiv} we show a relation between two types of infinite 
planar triangulations that demonstrates the universality principle. Through
this relation we also get an infinite form of the main result of \cite{BRW}
(see also \cite{BFSS}).

In a forthcoming paper (by the first author) \cite{UIPT2},
an alternative method of constructing
and sampling the UIPT is given. Using this method, it is shown there that
up to polylogarithmic factors the UIPT has growth rate $r^4$, agreeing both
with the heuristics for the Hausdorff dimension \cite{AmWa} and with the
asymptotics for the radius of finite maps \cite{ChSc}. That paper also
proves that the component of the boundary of the ball of radius $r$
separating it from infinity has size roughly $r^2$. The method also enables
an analysis of site percolation on the UIPT.

We proceed now to give formal definitions of the types of triangulations we
study. An exact formulation of our main results will follow.

%%%%%%%%%%%%%%%%%%%%%%%%
\subsection{Definitions}
%%%%%%%%%%%%%%%%%%%%%%%%

The notion of a triangulation is very similar to the topological notion of
a simplicial complex, although since we deal with the combinatorial aspects
rather then the topological ones we will use a graph theoretic approach.

The notion of a triangulation has a bit of ambiguity around it. There are
several variations on the definition, and they have much in common although
there are some minor differences between them. The common thread to all
variations is that a triangulation is a graph embedded in the sphere $S^2$
so that all faces are triangles. We will work with two types of
triangulations.

\begin{defn}\label{def:complex}
Consider a finite connected graph $G$ embedded in the sphere $S^2$.
A {\em face} is a
connected component of $S^2\setminus G$. The face is a {\em triangle} if
its boundary meets precisely three edges of the graph. Similarly, a face is
an {\em $m$-gon} if it meets $m$ edges. An 
{\em embedded triangulation} $T$ is such a
graph $G$ together with a subset of the triangular faces of $G$.

Let the {\em support} $S(T)\subset S^2$ of $T$ be the union
of $G$ and the triangles in $T$.
Two embedded triangulations $T,T'$ are considered {\em equivalent} 
if there is a 
homeomorphism of $S(T)$ and $S(T')$ that corresponds $T$ and $T'$.
$T$ is a triangulation of the sphere if $S(T)=S^2$.  It is
a triangulation of an $m$-gon if $S^2\setminus S(T)$ is
a single $m$-gon.
\end{defn}

For convenience, we usually abbreviate ``equivalence class of embedded
triangulations'' to ``triangulation''.  This should not cause much
confusion. The definition extends naturally to other manifolds, though we
will not be concerned with that generality here.

Following the terminology found in \cite{ADJ} for types of triangulations,
we define three classes of triangulations, types I, II and III.
These differ according to which graphs are permitted in
\ref{def:complex}. 
In type I, there may be more than one edge
connecting a pair of vertices, and loops (i.e., edges with
both endpoints attached to the same vertex) are allowed as well.
Type I triangulations will not be considered here,
though some of the results (and proofs) apply to them as well.

\begin{defn}\label{def:type2}
A {\em type II triangulation} is a triangulation where the underlying graph 
has no loops, but may have multiple edges.
\end{defn}

\begin{defn}\label{def:type3}
A {\em type III triangulation} is a triangulation where the underlying
graph is a simple graph (having no multiple edges or loops).
\end{defn}

Type II (resp.\ type III) triangulations are also referred to as 2-connected
(resp.\ 3-connected) triangulations, since they are the triangulations with
2 or 3 connected underlying graphs.

If $T$ is a
triangulation of a domain in the plane which may have several holes
(i.e., several boundary components), we will refer to the holes of the
domain as {\em external} or {\em outer} faces of $T$. An external face may
have 3 vertices on its boundary and then it is a triangle in itself. In
that case that face is
still distinguished from the triangles of $T$. In the case of type II, an
external face can also have only 2 vertices on its boundary.

It is worthwhile noting that the circle packing theorem \cite{Koebe} gives a
canonical embedding in the sphere (up to Moebius transformations) of a
type III triangulation of the sphere.

The vertices of $T$ lying on the boundary of its support $S(T)$
are called boundary vertices,
and those in the interior of $S(T)$ are internal vertices.
When we consider triangulations of a domain in the sphere with a number of
boundary components we will usually fix the number of boundary vertices in
each component as part of the domain. Thus, for example, a disc with $m$
boundary vertices will be distinguished from a disc with $m` \neq m$
boundary vertices. Such a disc is referred to as an $m$-gon.

The size of a triangulation $T$, denoted $|T|$, is defined as the number of
internal vertices. Since all faces are triangles, by Euler's characteristic
formula, if $E$ (resp.\ $F$) is the number of edges (resp.\ faces) of $T$, 
then $3|T|-E$ (resp.\ $2|T|-F$) is determined by the number and size of the 
boundary components of $|T|$. In particular, for a sphere all vertices are
internal, and so $3|T|-E=6$ and $2|T|-F=4$.

\begin{figure}
\begin{center}
\includegraphics[width=4.7in]{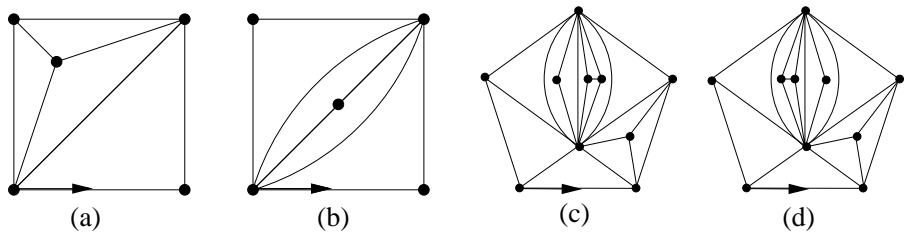}
\end{center}
\caption{\label{fig:eg1} A type III and three type II triangulations.
Triangulations (a) and (b) are triangulations of a square,
while (c) and (d) are triangulations of a pentagon.}
\end{figure}

Note that for a type III triangulation of the sphere (and even slightly
more generally) the underlying graph determines the triangulation,
i.e., whether any three edges form a triangle or not. When multiple edges
are allowed there may be several distinct embeddings of the graph in the
sphere giving distinct triangulations. E.g., in Figure~\ref{fig:eg1} (c)
and (d) are distinct triangulations that have the same underlying graph.

\bigskip

A fundamental problem encountered when studying planar maps (triangulations
included) is that of symmetries, namely that some maps have non trivial
automorphism groups. It seems plausible that most triangulations are
asymmetric. While this has been proved \cite{Tut2,RiWo2}, we dispose of
this problem in another manner. A simple way of eliminating any symmetries
there are is by adding a root to the triangulation.

\begin{defn}\label{def:root}
A {\em root} in a triangulation $T$ consists of a triangle $t$ of $T$ called
the {\em root face}, with an ordering of its vertices $(x,y,z)$. The vertex
$x$ is the {\em root vertex} and the directed edge $(x,y)$ is the
{\em root edge}.
\end{defn}

Note that in type II triangulations there may be more then one triangle
with the same three vertices, so marking only the three vertices does not
generally suffice. In a triangulation of the sphere, if the root edge is
given, then there are exactly two possibilities for the root. We will
usually mark only the root edge as in Figure~\ref{fig:eg1}, by an arrow.

When $T$ has a boundary we will usually assume that the root edge lies on
the boundary. Since a disc with $m$ boundary vertices is referred to as an
$m$-gon, triangulations (a) and (b) of Figure~\ref{fig:eg1} are of a
square, while (c) and (d) are of a pentagon.

\medskip

There are many possible variation on the definition of
triangulation. Restricting to 2 or 3-connected underlying graphs (or even 4
or 5-connected) gives slightly different definitions. It is possible to
restrict the degrees of vertices, or to allow faces that are not
triangles. Thus quadrangulations as well as convex polytopes in general,
can be described as a variation on the notion of triangulation. Most of the
results proved below should have analogues for such generalizations, though
the proofs do not always carry through and some complications are incurred
in the transition. For convenience and brevity, we will deal with type II
and type III triangulations here.

\medskip

The definition of an embedded infinite triangulations is identical to that of
finite triangulations, except that in that case $G$ is infinite, of course.
However, we will generally require of our infinite triangulations to
be locally finite, in the following sense.

\begin{defn}\label{def:loc_fin}
A triangulation $T$ embedded in the sphere $S^2$ is {\em locally finite} if
every point in $S(T)$
 has a neighborhood in $S^2$ that intersects only a finite number of
elements of $T$ (i.e., edges, vertices, triangles).
\end{defn}

We will henceforth require all the triangulations under discussion to be
locally finite.
In particular this implies that the graph of $T$ is locally finite (each
vertex is incident only to finitely many edges). However, this also
requires that the embedding of $G$ is faithful to the combinatorial
structure in the following sense: if $\{p_n:n=1,2,\dots\}$ is a sequence in
$S^2$ belonging to distinct edges of $G$, then accumulation points of $p_n$
must be outside of $S(T)$. The condition on the embedding is needed only so
that the definition of equivalence of embeddings, as stated in the finite
case, will also be useful in the infinite setting. It is possible
to define triangulations from a completely combinatorial standpoint, and
then embeddings may be ignored. However, the definitions are tricky, and it
is convenient to think of triangulations as embedded in the plane or sphere.

A triangulation of the plane is a triangulation $T$ with $S(T)=S^2\setminus
\{x\}$ for some $x$. The sphere may be identified with the plane with $x$
mapped to a point at infinity. Thus this definition is equivalent to having
$S(T)=\R^2$ with no accumulation points in $\R^2$.

\medskip

A triangulation may be endowed with a metric in a number of ways.
We will rather use a metric on the vertices of a
triangulation --- the graph metric induced by the underlying graph.
It is also interesting to consider a triangulation as a metric space
by having each face be isometric to an equilateral triangle with the
shortest path metric on the whole triangulation. Then a triangulation of
the sphere is a metric space homeomorphic to the sphere.

For either type, the space $\T$ of finite and infinite (equivalence classes
of) connected planar rooted triangulations is endowed with a natural
topology described in \cite{BeSc}. (Note however, that there the root was
only a vertex, leading to a closely related but slightly different notion).
A sequence of rooted triangulations converges to a triangulation $T$ if
eventually they are equivalent with $T$ on arbitrarily large combinatorial
balls around the root. This is a metric topology: e.g., set $d(T,T')=k^{-1}$,
where $k$ is the maximal radius such that the combinatorial balls of radius
$k$ around the roots are equivalent. In this topology, all finite
triangulations are isolated points, and infinite triangulations are their
accumulation points. It is not hard to check that $d$ is a metric; i.e.,
if $d(T,T')=0$, then $T$ is equivalent to $T'$ (even when $T$ and
$T'$ are infinite).

Given a Cauchy sequence of locally finite embedded rooted triangulations it is
easy to see that it
is possible to choose for them embeddings that eventually agree on the ball
of any fixed radius about the root.
Thus, the limit of the sequence exists
(as a locally finite embedded triangulation).
In other words, the space $\T$ of (locally finite embedded rooted)
triangulations is complete.

This metric is non-Archimedean, i.e.\ $d(T_1,T_2)\le\max\{d(T_1,T_3), 
d(T_2,T_3)\}$ which implies that if two balls in $(\T,d)$ intersect, then
one is a subset of the other.

Unlike in the setup of \cite{BeSc}, the triangulation space is not compact.
Consider the sequence $T_n$ of triangulations where $T_n$ contains two
vertices of degree $n$ with the same $n$ neighbors forming a cycle
(i.e., a double pyramid). Since $T_n$ are distinct and all have diameter 2,
$\{T_N\}$ has no convergent subsequence.

\medskip

We will be interested in the uniform distributions on triangulations:
\begin{defn}\label{def:taun}
$\tau^2_n$ (resp.\ $\tau^3_n$) is the uniform distribution on rooted type
II (resp.\ type III) triangulations of the sphere of size $n$ (i.e., having
$n$ vertices).
\end{defn}

The topology on the triangulation space induces a weak topology on the
linear space of measures supported on planar triangulations. We study the
distribution on infinite planar triangulations which is the weak limit of
$\tau_n$ as $n \to \infty$:

\begin{defn}\label{def:limit}
A measure $\tau$ on $\T$ is the limit of $\tau_n$ if for every 
bounded continuous
function $f:\T\to\R$
\[
\lim_{n \to \infty} \int f\, d\tau_n = \int f\, d \tau .
\]
\end{defn}

Since for any radius $r$ and triangulation $T$ the characteristic function
of the event $\{B_r=T\}$ is continuous, $\tau_n\to\tau$ implies convergence
with respect to neighborhoods of the root; that is, for any $r$ and $T$
\begin{equation}\label{eq:limit}
\tau_n\left(B_r=T\right) \to \tau\left(B_r=T\right) .
\end{equation}
However, since $\T$ is not compact, the existence
of $\lim_{n\to\infty} \tau_n(B_r=T)$ for every $r$ and $T$ is not sufficient
for convergence.
Indeed, if $\mu_n$ are distributions on triangulations where the
degree of the root vertex in $\mu_n$ is a.s.\ $n$,
then the probability of observing any
given triangulation as the ball of radius $r\ge 1$ tends to 0,
but the weak limit does not exist in the
sense of Definition~\ref{def:limit}.
(It would not even be correct to say that $\mu_n\to 0$,
since $\mu_n(\T)=1$.)

While the existence of the limits $\lim_n\tau_n(B_r=T)$
is not in itself sufficient for existence of a weak
limit of $\tau_n$, it follows from the equivalence
of Definition~\ref{def:limit} to convergence with respect to
the Prohorov metric (see \cite[\S11.3]{Dud}) that if
$\tau$ is a probability measure satisfying
(\ref{eq:limit}) for every $r$ and $T$,
then $\tau_n\to\tau$.
Indeed, for $\ep>0$
and $A\subset\T$ let $A^\ep$ denote the set of all $T\in \T$
such that $d(T,T')<\ep$ for some $T'\in A$.  Then $A^{1/r}$ is
just the set of triangulations $T$ such that the ball $B_r$ in
$T$ is the same as the ball $B_r$ in some $T'\in A$.
Given $r\in\{1,2,\dots\}$, we may find distinct $T_1,\dots,T_m$
such that $\tau\bigl(B_r\notin\{T_1,\dots,T_m\}\bigr)\le (2r)^{-1}$.
For all sufficiently large $n$ and for $j=1,2,\dots,m$, we have
$\tau(B_r=T_j)-\tau_n(B_r=T_j)<(2rm)^{-1}$.
Given a Borel set $A\subset\T$, let $J_A$ be the set
of $j\in\{1,\dots,m\}$ such that $A\cap \{B_r=T_j\}\ne\emptyset$.
Then
\begin{eqnarray*}
\tau(A)
& \le &
\tau(A^{1/r})\le
(2r)^{-1}+\sum_{j\in J_A} \tau(B_r=T_j)
\\
& \le &
(2r)^{-1}+\sum_{j\in J_A}\bigl( \tau_n(B_r=T_j)+(2 r m)^{-1}\bigr)
 \le
 \tau_n(A^{1/r})+r^{-1}
\end{eqnarray*}
 holds for all sufficiently
large $n$ and for all Borel sets $A\subset\T$.
Consequently, $\tau_n\to\tau$ in the Prohorov metric,
and therefore also in the sense of Definition~\ref{def:limit}.
Thus, understanding the limiting probabilities of balls is one of the key
ingredients in proving that $\lim_n\tau_n$ exists. 

\medskip

If $T_1$ and $T_2$ are rooted triangulations, we say that $T_1$ is
contained in $T_2$ (and write $T_1\subset T_2$) if the two roots are the
same and $T_1$ is contained in $T_2$ as unrooted triangulations.
Sometimes we may also write $T_1\subset T_2$ to mean that there is
a triangulation isomorphic to $T_1$ contained in $T_2$.

Finally, a word on notation.
By $X_n \sim Y_n$ we mean that $X_n/Y_n \to 1$.
By $X_n \approx Y_n$ we mean that $\log X_n/\log Y_n \to 1$.
We use $c,c_1,c_2,\dots$ to signify constants, whose actual value
may change from one formula to another.

%%%%%%%%%%%%%%%%%%%%%%%%%
\subsection{Main Results}
%%%%%%%%%%%%%%%%%%%%%%%%%

We will first prove that

\begin{thm}\label{thm:limit}
There exists a probability measure $\tau^2$ (resp.\ $\tau^3$) supported on
infinite planar triangulations of type II (resp.\ type III) such that
\[
\tau^i = \lim \tau^i_n .
\]
\end{thm}

\note The proof of this theorem basically has two parts.  The
first is to show that for every $r=1,2,\dots$ and every finite
triangulation $T$, the limits $\lim_{n\to\infty}\tau^i_n(B_r=T)$
exist.  This part is based on the exact enumeration results,
(and is not entirely new). Existence of these limits is a consequence of
the well behaved asymptotic forms for the enumeration results, and can be
extended to other classes of planar structures.

The second necessary ingredient is
to prove tightness (Lemma~\ref{lem:M_tight}), which is needed
since the space $\T$ is not compact.  To see why tightness is necessary,
the reader may wish to consider the sequence of probability measures
$\delta_n$ on $\N$ where
$\delta_n(A)=1_{n\in A}$.  This sequence does not converge to any measure.
However, $\lim_{n\to\infty}\delta_n(\{k\})=0$ for every $k\in\N$.
Similarly, it is easy to come up with examples of probability measures
$\mu_n$ on $\T$ where $\lim_{n}\mu_n(B_r=T)$ exists for every $T$
but the limit of $\mu_n$ does not exist.

\medskip

Having established the existence of the limit measure $\tau$
(we will often drop the type
notation when results hold for either type) we turn to study the a.s.\
properties of a sample of $\tau$. Denote such a sample by UIPT. A basic
geometric property, one endedness, will show that the limit structure 
maintains the plane's topology. Recall the definition:

\begin{defn}
A graph $G$ is said to have {\em one end} (is one-ended) if for any finite
subgraph $H$, $G \setminus H$ contains exactly one infinite connected
component.
\end{defn}

\begin{thm}\label{thm:one_end}
The UIPT is a.s.\ one ended, and is therefore a triangulation of
the plane.
\end{thm}

We also ask about the electrical type of the underlying graph.
In \cite{BeSc} it is shown that for any sequence of distributions on planar
graphs with degrees uniformly bounded by $M$, if a root is marked uniformly
in each graph then every subsequential limit is a.s.\ recurrent. This holds,
for example, for planar triangulations with uniformly bounded degrees.
However, for those distributions it is not clear how to prove that the
limit exists (simulations support this \cite{BCT}). The
following conjectures appear in \cite{BeSc}:

\begin{conj} \label{conj:limit}
For every $M\ge 6$, the distributions $\tau_N$ conditioned to have degrees
uniformly bounded by $M$ are weakly convergent.
\end{conj}

\begin{conj} \label{conj:type}
The UIPT is a.s.\ recurrent.
\end{conj}

VEL parabolicity (for vertex extremal length) is a property of infinite
graphs, closely related to circle packings for planar graphs. In graphs
with bounded degrees it is equivalent to recurrence \cite{HeSc}.
The proof in \cite{BeSc} of a.s.\ VEL parabolicity for
uniform infinite triangulations with bounded degrees is still valid for the
UIPT, with tightness (Lemma~\ref{lem:M_tight}) filling the role of bounded
degrees.

\medskip

\noindent{\bf Acknowledgments:} We thank Itai Benjamini and B\'alint
Vir\'ag for inspiring conversations. Part of this research was done during
visits of the first author to Microsoft Research. The first author thanks
his hosts for these visits.

\section{Counting} \label{sec:count}        %\input counting.tex

\subsection{Classical Results}

Much of the analysis of triangulations is based on counting them. This is
true both for finite triangulations and for infinite triangulations where
the asymptotics of the finite triangulations come into play.
The following counting
results go back to Tutte \cite{census} who counted various types of planar
maps and triangulations. The results we use here are not due to Tutte but
are derived using the same technique he uses. More details can be found in
\cite{Bro}. A good account of the technique including all results given
here can be found in \cite{GFBook}.

\begin{thm}\label{thm:count}
\be
\item For $n,m\geq 0$, not both 0, the number of type II 
triangulations of a disc with $m+2$ boundary vertices and $n$ internal
vertices that are rooted on a boundary edge is
\[
\phi^2_{n,m} = \frac {2^{n+1} (2m+1)! (2m+3n)!} {m!^2 n! (2m+2n+2)! } .
\]

\item For $m\geq 1, n\geq 0$, the number of rooted type III triangulations
of a disc with $m+2$ boundary vertices and $n$ internal vertices 
that are rooted on a boundary edge is
\[
\phi^3_{n,m} = \frac{2(2m+1)!(4n+2m-1)!}{(m-1)!(m+1)!n!(3n+2m+1)!} .
\]
\ee
\end{thm}

The case $n=m=0$ for type II triangulations warrants special attention. A
triangulation of a 2-gon must have at least one internal vertex so there
are no triangulations with $n=m=0$, yet the above formula gives
$\phi^2_{0,0}=1$. It will be convenient to use this value rather then 0 for 
the following reason. Typically, a triangulation of an $m$-gon is used not
in itself but is used to close an external face of size $m$ of some other
triangulation by ``gluing'' the boundaries together. When the external face
is a 2-gon, there 
is a further possibility of closing the hole by gluing the two edges to
each other with no additional vertices. Setting $\phi^2_{0,0}=1$ takes this 
possibility into account.

Since we will consider the asymptotics of large triangulations we will need 
the following estimates of these numbers. Using the Stirling formula, as
$n \to \infty$ we have the following:
\[
\phi^2_{n,m} \sim C^2_m \alpha_2^n n^{-5/2} ,
\]
where $\alpha_2=27/2$ and
\[
C^2_m = \frac{\sqrt{3}(2m+1)!}{4\sqrt{\pi}m!^2} (9/4)^m
        \sim C 9^m m^{1/2} .
\]

For type III triangulations we have similar estimates:
\[
\phi^3_{n,m} \sim C^3_m \alpha_3^n n^{-5/2} ,
\]
where $\alpha_3=256/27$ and as $m \to \infty$:
\[
C^3_m = \frac{2(2m+1)!}{6\sqrt{6\pi}(m-1)!(m+1)!} (16/9)^m
      \sim C (64/9)^m m^{1/2} .
\]

Much of the time we will not distinguish between type II and type III
triangulations. The type index will be dropped either when the stated
results hold for both types or when it is be clear which type is discussed.

We are interested in triangulations of the sphere that have no predefined
boundary. The number of those is given by:

\begin{prop}\label{prop:sphere}
For either type, the number of rooted triangulations of the sphere with $n$ 
vertices is $\phi_{n-3,1}$.
\end{prop}

\begin{proof}
Adding a triangle that closes
the outer face of a triangulation of a triangle makes a triangulation of
the sphere. Alternatively, removing the triangle incident on the root edge
that is not the root triangle gives a triangulation of a triangle rooted on
the boundary.
Thus, there is a bijection between triangulations of the sphere
with $n$ vertices and triangulations of a triangle with $n-3$ internal
vertices.
\end{proof}

We will also be interested in triangulations of discs where the number of
internal vertices is not prescribed. The following measure is of particular
interest:

\begin{defn} \label{def:free}
The {\em free distribution} on rooted triangulations of an $(m+2)$-gon,
denoted 
$\mu_m$, is  the probability measure that assigns weight 
\[
 \alpha^{-n} / Z_m(\alpha^{-1})
\]
to each rooted triangulation of the $(m+2)$-gon having $n$ internal vertices,
where
\[
Z_m(t) = \sum_n \phi_{n,m} t^n .
\]
\end{defn}

As before, $\mu^2_m$ (resp.\ $\mu^3_m$) will denote free type II (resp.\
type III) 
triangulations, and similarly for the partition functions $Z^2_m$ and
$Z^3_m$. Thus, the probability of a triangulation $T$, is proportional to
$\alpha^{-|T|}$, and $Z_m$ acts as a normalizing factor.

Note that by the asymptotics of $\phi$ as $n \to \infty$ we see that the
sum defining $Z$ converges for any $t \leq \alpha^{-1}$ and for no larger
$t$. The value of the partition functions will be useful. For this we have:

\begin{prop}\label{prop:Z}
\be
\item For type II triangulations, if $t=\theta(1-2\theta)^2$: 
\[
Z^2_m(t) = \frac{(2m)!((1-6\theta) m +2-6\theta)}{m!(m+2)!}
       (1-2\theta)^{-(2m+2)} .
\]

\item For type III triangulations, if $t=\theta(1-\theta)^3$: 
\[
Z^3_m(t) = \frac{(2m)!((1-4\theta)m+6\theta)}{m!(m+2)!}
       (1-\theta)^{-(2m+1)} . 
\]
\ee
\end{prop}

At the critical point $t=\alpha^{-1}$ we will omit $t$. There $Z$ takes
the values:
\[
Z^2_m = Z^2_m(2/27) = \frac{(2m)!}{m!(m+2)!} \left(\frac94\right)^{m+1} ,
\]
and
\[
Z^3_m = Z^3_m(27/256) = \frac{2(2m)!}{m!(m+2)!} \left(\frac{16}9\right)^m.
\]

The proof can be found as intermediate steps in the derivation of
$\phi_{m,n}$ in \cite{GFBook}. The above form may be deduced after a
suitable reparametrization of the form given there.

%%%%%%%%%%%%%%%%%%%%%%%%%
\subsection{Universality}
%%%%%%%%%%%%%%%%%%%%%%%%%

While the exponential term in the asymptotics of $\phi$ is different for
type II and III, the next term of $n^{-5/2}$ is the same. Similarities also
occur in the asymptotics of $C_m$ and of $Z_m$ for the two types. Those
similarities are not coincidental. It turns out that the asymptotic form is
quite common when counting 2 dimensional structures. That form of the
asymptotics is not dependent on the manifold, and is valid for any
2-dimensional manifold with or without boundaries. The same forms also
appear when instead of triangulations other types of maps are considered,
and was found to hold for a large variety of map types (\cite{BCT,BoKa} and
also the result of \cite{ChSc}, related to our growth results). We therefore
believe that many of the results on the UIPT hold in a much more general
context. In this, infinite planar objects are similar to random walks,
critical percolation and many other critical models where the large scale
properties are independent of the local lattice.

This universality is related to the basic property of the $2$-sphere that a
cycle partitions it into two parts, i.e., the Jordan Curve Theorem. This
leads to a similarity 
between recurrence relations for different types of structures and through
them to similar asymptotics for the solutions. For another instance of
universality and some explanation see \cite{BFSS}.

It turns out that the exponential part of the asymptotics will cancel out
often and when finer properties of infinite triangulations are considered
the power term will come into play and determine the observed behavior.

\subsection{Some Estimates}

We will need the following estimates throughout the paper.

\begin{lemma}\label{lem:small_sum}
Let
\[
S(k,N,a) = \sum_{\substack{n_1+\ldots+n_k=N \\ n_1,n_2 > a}}
     \left( \prod n_i \right)^{-5/2} ,
\]
then for any $k$ there is a $c=c(k)$ such that for any $N$:
\[
S(k,N,a) \leq c N^{-5/2} a^{-3/2} .
\]
\end{lemma}

\begin{proof}
Clearly
\[
S(k,N,a) \leq k! \sum_{\substack{n_1\geq \ldots \geq n_k \\
                    n_1+\ldots +n_k=N \\ n_1,n_2 > a}}
     \left( \prod n_i \right)^{-5/2} ,
\]
since each term in the sum over ordered $k$-tuples corresponds to at most
$k!$ terms in the original sum, and less if there are any repetitions. 
Since each possible choice of $n_2,\ldots,n_k$ determines a unique value
for $n_1$ and always $n_1 \geq N/k$ we can replace $n_1$ by the smaller
$N/k$ and extend the range of summation. 
\begin{eqnarray*}
S(k,N,a)
    & \leq & k! (N/k)^{-5/2}
            \sum_{\substack{n_2\geq \ldots \geq n_k \\ n_2 > a}}
                    \left( \prod_{i\neq 1} n_i \right)^{-5/2}           \\
    & \leq & k! (N/k)^{-5/2} \left( \sum_{n_2\geq a} n_2^{-5/2} \right)
            \prod_{i>2} \left( \sum_{n_i} n_i^{-5/2} \right)        \\
    & \leq & c N^{-5/2} a^{-3/2} .
\end{eqnarray*}
\end{proof}

\section{Basic Properties} \label{sec:basic}    %\input basic.tex

%%%%%%%%%%%%%%%%%%%%%%%%%%%%%%%%%%%%%%%%%%%
\subsection{Invariance with respect to the Random Walk}
%%%%%%%%%%%%%%%%%%%%%%%%%%%%%%%%%%%%%%%%%%%

If we are given a finite triangulation, but not the location of
the root, what can we say about the location of the root? The following
proposition says that not much. For a triangulation $T$ and a possible root
$r$ in $T$ let $T_r$ denote the triangulation $T$ with $r$ marked as root
(if $T$ is rooted then the old root is no longer marked).

\begin{prop} \label{prop:reroot}
Let $T$ be a sphere triangulation chosen by $\tau_n$, and $r$ be a root in $T$ 
chosen uniformly among all possible roots. Then $T_r$ is uniformly
distributed among all rooted triangulations (of size $n$). 
\end{prop}

\begin{proof}
At first glance this seems trivial: since all rooted triangulations are
equally likely no triangle in $T$ should be more likely to be the root than
any other. However, there is a subtlety here since there may be several
triangles $r$ such that the triangulations $T_r$ are isomorphic. This
occurs whenever $T$ has a non trivial automorphism.

The key fact here is that any automorphism of $T$ that preserves a root is
necessarily the identity automorphism. If $R$ is the set of possible roots
and $G$ is the automorphism group of $T$ as an unrooted triangulation, then
$G$ acts naturally on $R$ and a non identity element of $G$ has no fixed
points in $R$.
Thus, the size of the orbit of a triangle $r \in R$ is just the size
of $G$, regardless of $r$.

Since each of the orbits in the action of $G$ on $R$ corresponds to a
distinct rooted triangulation, and each orbit has the same size, each
possible triangulation is equally likely to result after a new root is
selected. 
\end{proof}

Note that since each directed edge can be completed in two ways to a root
each directed edge is equally likely to be the root edge. From this we see
that the UIPT must be invariant with respect to a random walk:

\begin{thm}\label{thm:RW}
Let $T$ be a triangulation chosen by $\tau_N$ for some $N$ or by a
subsequential limit $\tau$. If $x$ is the root vertex of $T$, $y$ is a
uniformly chosen neighbor of $x$, and $(y,z,w)$ is a triangle in $T$
uniformly chosen among all triangles including $y$, then $T_{(y,z,w)}$ has
the same law as $T$.
\end{thm}

\begin{proof}
For finite $N$, if a vertex $x$ of degree $d$ is the root vertex, then
there are $d$ possibilities for the root edge (and $2d$ options for the
root). It follows that the probability that $x$ is the root is proportional
to its degree. This is the stable distribution for the random walk on the
graph of $T$, so as a consequence of Proposition~\ref{prop:reroot} we see
that $T_{(y,z,w)}$ has the same law as $T$.

Since this is true for every $\tau_N$, the same holds for any subsequential
limit.
\end{proof}

%%%%%%%%%%%%%%%%%%%%%%%%%%
\subsection{One Endedness}
%%%%%%%%%%%%%%%%%%%%%%%%%%

We start with a lemma describing the behavior of a triangulation on
a disjoint union of discs.

\begin{lemma} \label{lem:one_big}
Given $k$ disjoint polygons (with given boundary sizes)
and a triangulation $T$ of the
polygons, let $n_i$ be the number of internal vertices in the $i$'th polygon.
Then
\[
\left|\bigl\{T\bigm|\sum n_i=N \wedge \exists i,j,\ i\ne j,\
 n_i,n_j>a \bigr\}\right|
 < C \alpha^N N^{-5/2} a^{-3/2} ,
\]
where $C$ depends only on the number and sizes of the boundaries
of the polygons.
\end{lemma}

\begin{proof}
We prove that the number of triangulations where $n_1,n_2>a$ is small, as
required. By symmetry, the number for any other pair $(i,j)$ has the same
bound. Since the number of such pairs, $\binom k 2$, does not depend on $a$
or on $N$, this suffices.

We use the upper bound $\phi_{n,m} \leq \beta_m (n+1)^{-5/2} \alpha^n$
($+1$ is only necessary to account for $n=0$, and is not essential). Assume 
the $i$'th domain has boundary size $m_i+2$. The number of triangulations
we wish to bound is:
\begin{eqnarray*}
  \sum_{\substack{n_1 + \ldots + n_k = N \\ n_1,n_2 > a}}
    \prod_i \phi_{n_i,m_i}
& \leq &
  \sum_{\substack{n_1 + \ldots + n_k = N \\ n_1,n_2 > a}}
    \prod_i \beta_{m_i} (n_i+1)^{-5/2} \alpha^{n_i}     \\
& = &
  \alpha^N \prod \beta_{m_i}
  \sum_{\substack{n_1 + \ldots + n_k = N \\ n_1,n_2 > a}}
    \prod (n_i+1)^{-5/2}                    \\
& \leq &
  c_1 \alpha^N
  \sum_{\substack{n_1 + \ldots + n_k = N+k \\ n_1,n_2 > a}}
    \prod (n_i)^{-5/2}                      \\
& \leq &
  c_2 \alpha^N N^{-5/2} a^{-3/2} ,
\end{eqnarray*}
where at the end we used Lemma~\ref{lem:small_sum}.
\end{proof}

Generally a limit $T$ of a sequence of finite sphere triangulations need
not have support $S(T)$ which is homeomorphic to the sphere or even the plane.
While the limit is still
planar, when embedded in the sphere $S(T)$ may have any number of accumulation
points. One accumulation point gives a punctured sphere, i.e., the
plane. More than one means that $S(T)$ has a more complicated
topological structure; it is no longer simply-connected.

\begin{coro} \label{cor:one_end}
Every subsequential limit of $\tau_N$ a.s.\ has one end.
\end{coro}

\begin{proof}
Suppose that a subsequential limit $\tau$ has more then one end with
positive probability. Then for some $k$ and some $\ep>0$ the probability
that a loop of length $k$ including the root partitions a sample of $\tau$ 
into two infinite parts is at least $\ep$. This implies that for any $a$
for infinitely many $N$ the $\tau_N$-probability of having a loop of length $k$
including the root that has at least $a$ vertices on either side is at
least $\ep/2$. Call such a loop a separating loop.

Count pairs $(T,L)$ with $T$ a triangulation of size $N$ and $L$ a
separating loop included in $T$. From Lemma~\ref{lem:one_big} we know that
the total number of such pairs is $O(\alpha^N N^{-5/2} a^{-3/2})$. However,
the total number of sphere triangulations with $N$ vertices is
$\phi_{N-3,1} \sim C \alpha^N N^{-5/2}$, and by dividing we deduce that the
expected number of separating loops is $O(a^{-3/2})$. In particular as
$a \to \infty$ so does the probability that a separating loop exists.
\end{proof}

\section{Existence of the Limit} \label{sec:limit}
\subsection{Tightness}              %\input tight.tex

Given the formula for the number of disc triangulations (or even its
asymptotics), it is simple to verify that the probability of observing any
given ball around the root converges to some limit. However, since $\T$ is
not compact this is not sufficient to guarantee convergence of the
measures. The missing factor is a tightness result for the measures.

Recall that a family of random variables $\{X_n\}$ is tight with respect to
$n$ if 
\[
\lim_{t \to \infty} \P(|X_n| > t) = 0
\]
uniformly with respect to $n$. More generally, a family of probability 
measures $\{\mu_n\}$ on a topological space is tight if for any $\ep>0$
there is a compact set $A$ with $\mu_n(A)>1-\ep$ for all $n$.

We first need to prove that the uniform measures $\{\tau_n\}$ on finite
triangulations are a tight family. To this end, we first prove the
following estimates for the degree of the root in either type of
triangulation. While the lemmas are very similar in nature, the methods of
proof given here are different. This demonstrates the underlying unity of
the different models, while local differences make some techniques
applicable in one and others in another. The following two Lemmas appear in
similar form in \cite{GaRi}. 

\begin{lemma} \label{lem:deg_dist1}
Denote the degree of the root vertex by $d_0$. For any $\ep>0$ there is a
$c=c(\ep)$ such that
\[
\tau^3_N(d_0=k) < c \left(\frac{3}{4}+\ep \right)^k,
\]
uniformly for all $N$ and
\[
\tau^3_N(d_0=k) \xrightarrow[N \to \infty]{unif}
  \frac{2(2k-3)!}{(k-3)!(k-1)!} \left( \frac{3}{16} \right)^{k-1} .
\]
\end{lemma}

\begin{proof}
A type III triangulation of the sphere where the root vertex
has degree $k$ is the
union of two triangulations, $T_0,T_1$ whose intersection is a $k$-gon:
$T_0$ contains the root
vertex and $k$ triangles connecting it to the sides of the
$k$-gon, and $T_1$ contains all other triangles. The root triangle has one
edge in the intersection of $T_0$ and $T_1$.
Choose this edge to be the root edge of
$T_1$. Now $T \leftrightarrow T_1$ is a bijection between
rooted triangulations of the sphere with $d_0=k$ and rooted triangulations
of a $k$-gon with the root edge on the boundary.

If $|T|=N$, then $|T_1|=N-k-1$, and we know
the number of such triangulations. Dividing by the number sphere
triangulations, we get:

\begin{eqnarray*}
\tau_N(d_0 = k)
 & = & \frac{\phi_{N-1-k,k-2}}{\phi_{N-3,1}}        \\
 & = & \frac{2(2k-3)!}{3(k-3)!(k-1)!}
       \frac{(4N-2k-9)!(3N-6)!(N-3)!}{(4N-11)!(3N-k-6)!(N-k-1)!}    \\
 &\to& \frac{2(2k-3)!}{(k-3)!(k-1)!} \left( \frac{3}{16} \right)^{k-1} .
\end{eqnarray*}

To prove the uniform exponential bound consider the ratio
\[
\frac{\tau_N(d_0=k+1)}{\tau_N(d_0=k)}
 = \frac{(2k-1)(2k-2)(3N-k-6)(N-k-1)}{k(k-2)(4N-2k-9)(4N-2k-10)}
 < \frac{3}{4}+\ep 
\]
for any sufficiently large $N,k$. 
\end{proof}

For type II, since multiple edges are present, there are two notions of
degree. The vertex degree of $v$ is the number of neighbors it has, while
the edge degree of $v$ is the number of edges incident on it. For our
purposes bounding the vertex degree is sufficient, but in what follows we
bound the larger edge degree.

\begin{lemma} \label{lem:deg_dist2}
Denote the edge degree of the root vertex by $d_0$, then for any $\ep>0$
there is a $c=c(\ep)$ such that
\[
\tau^2_N(d_0=k) < c \left(\frac{5}{3\sqrt 3} \right)^k 
\]
uniformly for all $N$.
\end{lemma}

\begin{proof}
Let $t_1,\ldots,t_{d_0}$ be the triangles incident with the root vertex,
ordered counterclockwise starting with the root triangle $t_1$. For
$s \leq d_0$ let $T_s$ be the sub-triangulation including triangles
$t_1,\ldots,t_s$. Adding $t_i$ one at a time, we consider the distribution
of $T_{s+1}$ conditioned on $T_s$, and show that for any $T_s$ there is a
probability bounded away from 0 that
$d_0=s+1$.

$T_s$ may have several external faces. One of those, say $F$, includes the
root vertex, and $t_{s+1}$ is in $F$. In order for $t_{s+1}$ to be the last
triangle adjacent to the root vertex it must include the two edges of $F$
on either side of the root vertex. Thus, in Figure~\ref{fig:tight}(a), the
triangles incident with the root vertex are numbered.
The triangle $t_{10}$ is the final
triangle, and it includes both the edge from $t_9$ and the edge from the
root triangle $t_1$. Note that when triangle $t_5$ is added, an unknown part of
the triangulation is enclosed, but this will not effect the bounds we get.

\begin{figure}
\begin{center}
\includegraphics[width=4.7in]{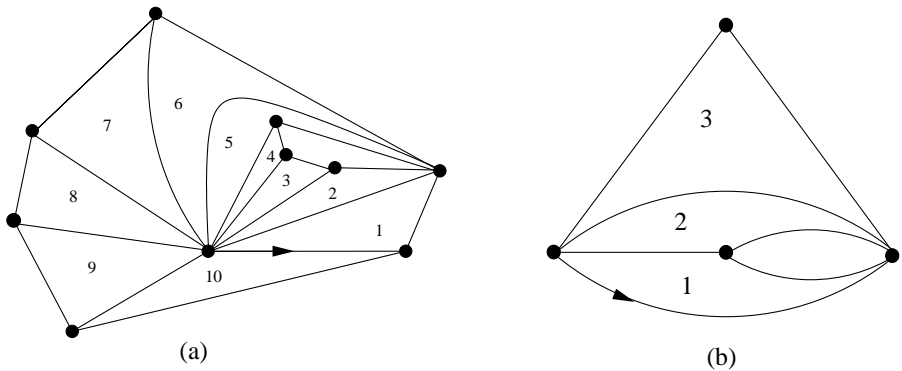}
\end{center}
\caption{\label{fig:tight} Proof of the tightness of the root's degree.}
\end{figure}

To bound from below
the probability that $t_{s+1}$ is the last triangle conditioned on 
$T_s$, assume that the boundary of $F$ has size $m+2$. At first condition on the
event that the part of the triangulation inside $F$ has $n$ vertices. The
number of possible ways to triangulate $F$ under these constraints is
$\phi_{n,m}$. If $t_{s+1}$ is the last triangle around the root, then adding
it leaves a face of boundary size $m+1$ with $n$ internal vertices. Thus
the probability of the next triangle being the last one is: 
\begin{eqnarray*}
\frac{\phi_{n,m-1}}{\phi_{n,m}}
  &=& \frac{m^2(2n+2m+1)(2n+2m+2)}{2m(2m+1)(3n+2m-1)(3n+2m)}    \\
  &>& \frac{2m(n+m)^2}{(2m+1)(3n+2m)^2} .
\end{eqnarray*}

If $m>0$, then this is at least 2/27, so the probability that $d_0>s+1$ is
at most 25/27. Since this bound is uniform it also holds when conditioning
only on $T_s$ and not on the number of internal vertices in $F$.

Thus, as new triangles are revealed, each triangle has a probability of at
least 25/27 of being the last one, unless $m=0$. If $m=0$, as for $T_2$ in
Figure~\ref{fig:tight}(b), then after a triangle is added we must have
$m=1$ and so out of every two consecutive $s$, at least one has $m>0$.
It follows that the probability of having more than $k$ edges leaving the
root vertex is at most $(25/27)^{(k-1)/2}$, as claimed.
\end{proof}

\note 
For type III triangulations, Lemma~\ref{lem:deg_dist1} gives the exact
probability of any given
degree in the UIPT. To a large extent, this is possible because the radius 1
neighborhood of the root has a simple structure. When multiple edges are
allowed, even the ball of radius 1 around the root can have a complicated
structure, making an exact calculation harder to get. On the other hand,
for type II triangulations, we can calculate the exact probability that a
certain triangle is present in the triangulation conditioned on some
sub-triangulation (e.g., the probability that $t_{s+1}$ is the last
triangle around the root conditioned on $T_s$, as in the proof). This is much
harder to do for type III triangulations, because we need to keep track of
which pairs of vertices already have edges between them, whereas in type II 
triangulations adding another edge is always legal.

\medskip

At this point we will rigorously define the ball $B_r$ of radius $r$ around
the root (or any other vertex, for that matter). This ball is a
sub-triangulation, but there is some subtlety in its definition. The
vertices of $B_r$ are all those vertices at distance at most $r$ from the
root vertex, but not all edges and triangles between these vertices are
necessarily part of $B_r$.

\begin{defn}
$B_0$ is just the root vertex itself. $B_{r+1}$ is composed of all
triangles incident on any vertex of $B_r$ together with their vertices and
edges.
\end{defn}

Note that there may be edges between vertices on the boundary of $B_r$ that
are not part of $B_r$ itself. Next, we turn our attention to the size of the
ball $B_r$. The following lemma holds for both types.

\begin{lemma} \label{lem:M_tight}
For any fixed $r$ the random variables $M_r = \max \{d_v \mid v \in B_r\}$
(i.e., the maximal degree in $B_r$) defined on $\T$ with measure $\tau_N$
are tight with respect to $N$.
\end{lemma}

\begin{proof}
For $r=0$, $B_r$ is just the root, and Lemmas~\ref{lem:deg_dist1} and
\ref{lem:deg_dist2} show that the degree of the root is tight with respect
to $N$ for either type.

We proceed by induction on $r$. Suppose that $M_r$ is tight with respect to
$N$. To show that $M_{r+1}$ is also tight we use Theorem~\ref{thm:RW}. Let
$T$ denote a sample of $\tau_n$, and let $X_0,X_1,\ldots$ be a simple random
walk on $T$ started at the root vertex $X_0$. Denote by $\P$ the resulting
probability measure on triangulations with paths beginning at the root. It
follows from Theorem~\ref{thm:RW} that for any $i$ the degree of $X_i$ has
the same distribution as the degree of the root. Fixing $M'>M>0$ we
estimate the probability that $M_r\le M$ and yet $M_{r+1}>M'$. Conditioned
on this event, there is at east one vertex $u\in B_{r+1}\setminus B_r$ with
$d_u>M'$. Since there is a path of length $r+1$ from the root vertex to
$u$, and all vertices on the path have degrees at most $M$,
\[
\P(d_{X_{r+1}}>M' \mid M_r\le M<M'<M_{r+1}) \ge M^{-(r+1)} ,
\]
and so
\[
\P(d_{X_{r+1}}>M') \ge M^{-(r+1)} \P(M_r\le M<M'<M_{r+1}) .
\]
By Theorem~\ref{thm:RW} the LHS does not depend on $r$ and is simply
$\tau_n(d_0>M')$. The RHS does not depend on the random walk either, so for
any $M$
\[
M^{r+1} \tau_n(d_0>M') \ge
\tau_n(M_{r+1}>M') - \tau_n(M_r > M) .
\]
By induction, for all $\ep>0$ we may choose $M=M(\ep)$
such that $\tau_n(M_r > M)\le\ep/2$ for all $n$.
Then we take $M'=M'(\ep)>M$ sufficiently large
so that $\tau_n(d_0>M')<M^{-(r+1)}\,\ep/2$
for all $n$.  This gives $\tau_n(M_{r+1}>M')<\ep$ for all $n$, and
completes the proof.
\end{proof}

\begin{coro} \label{cor:tight}
The family of probability measures $\{\tau_N\}$ is tight. 
\end{coro}

\begin{proof}
Since for any $r,L$ there are only finitely
many radius $r$ balls with $M_r<L$, if $L_r$ is any sequence,
then the set $A$ of
triangulations with $M_r<L_r$ is compact. From tightness of $M_r$ it follows
that for any $\ep$ there is a sequence $L_r$ so that $\tau_N(A) > 1-\ep$
for all $N$.
\end{proof}

Since any tight family of probability measures on a complete separable
space has a converging subsequence (see \cite[Theorem 11.5.4]{Dud}),
this implies 

\begin{coro} \label{cor:subseq}
Every subsequence of $\tau_N$ has a subsubsequence converging to a
probability measure. 
\end{coro}

\subsection{Taking the Limit}

In \cite{RiWo1} it is shown that for every finite triangulation $T$ there
is a constant $c$ such that asymptotically, in almost every sphere
triangulation of size $n$ the number of times $T$ appears is roughly $cn$.
This $c(T)$ is roughly the probability that a neighborhood of the root in
the UIPT is isomorphic to $T$. In fact, the result of \cite{RiWo1} is
stronger, since it gives not just an annealed probability of seeing $T$ but
that the quenched probability is constant. We bring here a simpler
calculation just for the annealed probability, since the results of the
calculation are useful in what follows.

It will be easier to work with rooted triangulations $A$ having the
property that if they are a sub-triangulation of the UIPT, then they appear
in it exactly once. This is not only the case: a root triangle together
with a cycle of some length may appear in the triangulation in several
different ways. 

\begin{defn} \label{def:rigid}
A rooted triangulation $A$
is {\em rigid} if it is connected and no triangulation includes two
distinct copies of $A$ with coinciding roots.
\end{defn}

The balls $B_r$ of a triangulation are rigid, as is evident from the
following sufficient criterion for rigidity (the proof is left to the
reader).

\begin{lemma}\label{lem:rigid}
If in the dual graph of triangulation $A$ the vertices corresponding to the
triangles form a connected set, and every vertex of $A$ is incident on a
triangle, then $A$ is rigid.
\end{lemma}

\begin{figure}
\begin{center}
\includegraphics[width=4.7in]{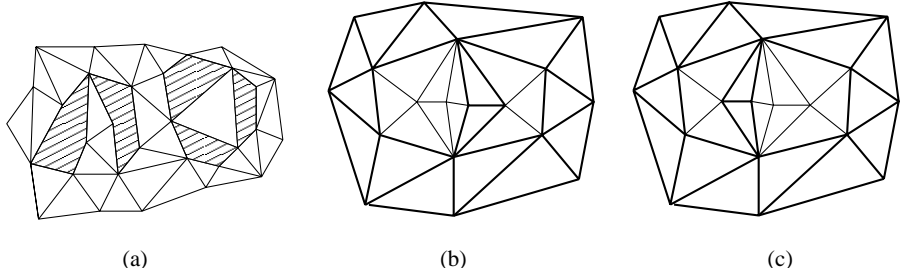}
\end{center}
\caption{\label{fig:rigid} A rigid triangulation (with shaded outer faces)
and two isomorphic completions of a non-rigid triangulation.}
\end{figure}

This criterion is not necessary for rigidity, as is demonstrated by
Figure~\ref{fig:rigid}(a), where there is an isolated triangle. In fact a
sufficient and necessary criterion is that the support $S(T)$
 be 3-connected. In order to
complete a planar triangulation to a sphere triangulation we need to fill
each of its external faces with some triangulation. The advantage of rigid
triangulations is that filling the external faces in different ways must
lead to distinct sphere triangulations, whereas for non-rigid
triangulations different ways of filling the faces may give rise to the
same complete triangulation. Figure~\ref{fig:rigid}(b,c) give an example of
a non-rigid triangulation and how two completions give rise to the same
triangulation.

A second consequence of the construction of $B_r$, whose proof (an
application of the Jordan curve theorem) is left to the reader is:

\begin{lemma}\label{lem:noedge}
In the ball $B_r$ there are no edges between two vertices of any external
face except those making the face itself.
\end{lemma}

\begin{prop}\label{prop:sub_prob}
Let $A$ be a rigid rooted triangulation having no edges between two
vertices of an external face except those making the face itself. Assume
$A$ has $n$ vertices, some of which are on $k$ external boundary components
of sizes $m_1+2,\ldots,m_k+2$. Then every subsequential limit $\tau$ of
$\tau_N$ has:
\[
\tau(A \subset T) = \frac{\alpha^{3-n}}{C_1} \Bigl(\prod_{i=1}^k Z_{m_i}\Bigr)
            \sum_{i=1}^k \frac{C_{m_i}}{Z_{m_i}} .
\]
Moreover, the probability that the $i$'th face is the infinite one
corresponds to the $i$'th term in the sum, i.e.:
\[ 
\frac{\alpha^{3-n}}{C_1} C_{m_i} \prod_{j \neq i} Z_{m_j} .
\]\end{prop}

\note For type II triangulations the restriction on edges between
vertices of an external face is not necessary. For type III triangulations
it is needed, since when such an edge exists it imposes restrictions on the
component inside the face. In general, the probability that $A \subset T$
can be found using the proposition together with the inclusion-exclusion
principle. The requested probability is a linear combination of a fixed
number of terms and each of them has a limit as above.

\begin{proof}
Let $\tau$ be a subsequential limit of $\tau_N$.
Denote by $Q=Q(A,n_2,\ldots,n_k)$ the event that $A \subset T$ (with the A's
root corresponding to $T$'s root) and that the part of $T$ in the $i$'th
external face of $A$ contains $n_i$ internal vertices. This is defined in
the finite as well is the infinite setting
(though we keep $n_2,\dots,n_k<\infty$). In what follows $n_1$ denotes
the number of vertices in the 1st external face, i.e.,
$n_1=N-n-\sum_{i>1} n_i$. The probability of $Q$ is:  
\[
\tau_N(Q) =
   \frac{\prod_{i=1}^k \phi_{n_i,m_i}}{\phi_{N-3,1}} ,
\]
and $\tau_N(A \subset T)$ is the sum over all possible vectors $(n_i)$ of
this probability.

We first consider the limit:
\begin{eqnarray*}
\lim_{N \to \infty} \tau_N(Q(A,n_2,\ldots,n_k))
  &=& \left(\prod_{i>1} \phi_{n_i,m_i}\right)
      \lim_{N \to \infty} \frac{\phi_{n_1,m_1}}{\phi_{N-3,1}}       \\
  &=& \left(\prod_{i>1} \phi_{n_i,m_i}\right)
      \lim_{N \to \infty} \frac{C_{m_1} n_1^{-5/2} \alpha^{n_1}}
                   {C_1 (N-3)^{-5/2} \alpha^{N-3}}      \\
  &=& \left(\prod_{i>1} \phi_{n_i,m_i}\right)
      \frac{C_{m_1}\alpha^{-n-\sum_{i>1} n_i}}{C_1 \alpha^{-3}} .
\end{eqnarray*}

Since the limit exists, it equals $\tau(Q)$.  This may be written as:
\begin{equation} \label{eq:Q_prob}
\tau(Q) = \frac{\alpha^{3-n}}{C_1} C_{m_1}
    \prod_{i>1} \phi_{n_i,m_i} \alpha^{-n_i} .
\end{equation}
Of course, a similar expression holds when the role of the 1st face is
filled by some other face, i.e., the sizes of all but the $i$'th face are
fixed.

Let $R_i=R_i(A)$ denote the event that $A\subset T$ and
all the external faces of $A$ except possibly the $i$'th one contain
finitely many vertices.  Obviously,
\[
R_1 = \bigcup_{n_2,\ldots,n_k<\infty} Q(A,n_2,\ldots,n_k) . 
\]
Using (\ref{eq:Q_prob}) we get for any subsequential limit $\tau$ of $\tau_N$:
\begin{equation} \label{eq:R_prob}
\tau(R_1) = \frac{\alpha^{3-n}}{C_1}C_{m_1} \prod_{i>1} Z_{m_i},
\end{equation}
and a similar formula for $R_j$, $j>1$.
It is clear that $\tau(R_i\cap R_j)=0$ for $i\ne j$ in $\{1,\dots,k\}$.
Moreover, Corollary~\ref{cor:one_end} (one end) implies 
$\tau\bigl(\{A\subset T\}\setminus\cup_i R_i\bigr)=0$.
Hence,
\[
\tau(A \subset T) = \sum_{i=1}^k \tau(R_i)
= \frac{\alpha^{3-n}}{C_1} \left(\prod Z_{m_i}\right)
            \sum \frac{C_{m_i}}{Z_{m_i}}\,.
\]
\end{proof}

\begin{proof}[Proof of Theorem~\ref{thm:limit}]
Corollary~\ref{cor:subseq} tells us that a subsequential limit $\tau$ exists.
It remains to prove that the limit does not depend on the subsequence.
Since the balls $B_r$ are rigid, Proposition~\ref{prop:sub_prob}
shows that for every $A$,   $\tau(B_r=A)$ does not depend on the subsequence.
Since the sets $\{T\mid B_r=A\}$ form a basis for the topology on
$\T$, this implies that $\tau$ does not depend on the subsequence.
\end{proof}

\begin{proof}[Proof of Theorem~\ref{thm:one_end}]
This follows from Theorem~\ref{thm:limit} and Corollary~\ref{cor:one_end}.
\end{proof}

\section{Locality}    \label{sec:local}     %\input locality.tex

Next, we look at another basic property of the UIPT, namely locality. The
meaning of locality is that isolated regions of the UIPT are almost
independent.  In the following, $R_i=R_i(A)$ will denote the event
defined in the proof of Proposition~\ref{prop:sub_prob}.

\begin{thm}\label{thm:locality}
Let $A$ be a finite rigid triangulation (for type III, with no edges
between vertices on external faces). Assume $A$ has $k$ external faces of
sizes $m_1+2,\ldots,m_k+2$. Condition on the event $R_i(A)$, and let
$T_j$ denote the component of the UIPT in the $j$'th face. Then:
\be
\item
The triangulations $T_j$ are independent.
\item
$T_i$ has the same law as the UIPT of an $(m_i+2)$-gon
(that is, the $N\to\infty$ limit of the uniform measure
on rooted triangulations of an $(m_i+2)$-gon with $N$ internal vertices).
\item
For $j \ne i$, $T_j$ has the same law as the free triangulation of an
$(m_j+2)$-gon. 
\ee
\end{thm}

\begin{proof}
Without loss of generality,
 assume $i=1$. From Equations~\ref{eq:R_prob} and \ref{eq:Q_prob} we
see that 
\[
\tau(R_1) = \frac{\alpha^{3-n}}{C_1} C_{m_1} \prod_{j>1} Z_{m_j}
\]
and so
\[
\tau\bigl(|T_j|=n_j \textrm{ for $j>1$}\bigm| R_1\bigr) =
    \prod_{j>1} \frac{\phi_{n_j,m_j} \alpha^{-n_j}}{Z_{m_j}} .
\]

Thus, we see that conditioned on $R_1(A)$ the sizes of the $T_j$'s
are independent, and $|T_j|$ is distributed like the free triangulation of
an $(m_j+2)$-gon. 
Consider $\tau_N$. Conditioned on
$Q(A,n_2,\ldots,n_k)$, since all possible triangulations of the sphere
with the prescribed component sizes are
equiprobable, the same holds for each component $T_i$. Thus, for any $N$ the
joint distribution of $(T_1,\ldots,T_j)$ conditioned on their sizes is a
product distribution. As $N \to \infty$, these joint distributions converge
to the product distribution, where $T_j$ is uniform on triangulations
with $|T_j|=n_j$.

Finally, the marginal of $T_1$ has size tending to infinity, and so
converges to the UIPT of an $(m_1+2)$-gon.
\end{proof}

\section{Ball Structure}   \label{sec:GW}       %\input GW_proc.tex

Recall that Theorem~\ref{thm:locality} tells us that conditioned on a
sub-triangulation $T$, with some external faces, the probability that a
face of size $m+2$ is the infinite one is proportional to
$\frac{C_m}{Z_m}$. In the case of type II or III triangulations we have:
\[
\frac{C^2_m}{Z^2_m} = \frac{(m+1)(m+2)(2m+1)}{3\sqrt{3\pi}} ,
\]\[
\frac{C^3_m}{Z^3_m} = \frac{m(m+2)(2m+1)}{6\sqrt{6\pi}},
\]
so in either case the probability of a face of size $m$ being the
infinite face is roughly proportional to $m^3$.

\begin{figure}
\begin{center}
\includegraphics[height=1.8in]{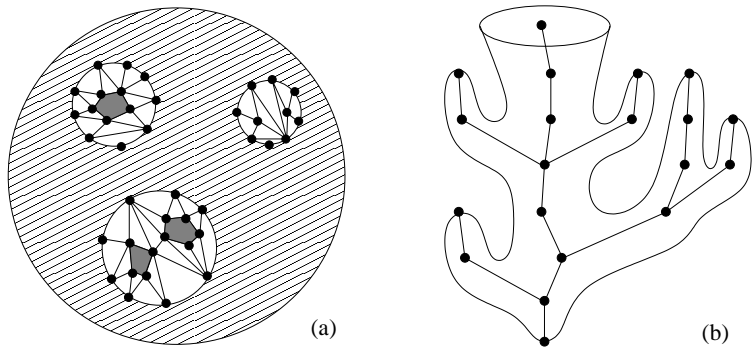}
\end{center}
\caption{\label{fig:ball} (a) A possible ball in a planar triangulation.
(b) A tree corresponding to a surface. Height corresponds to distance from
the root.}
\end{figure}

\medskip
We wish to study the relation between the ball of radius $r$ and the
ball of radius $r+1$. The ball of radius $r$ is a finite triangulation with 
any number of external faces with any combination of boundary sizes. Moving
to $r+1$ we add in each outer face some triangles around its circumference.
These added triangles can fill up the face, or they can split that face up
into a number of sub-faces of different sizes. Figure~\ref{fig:ball}(a)
shows a ball with several finite faces, and the layer of the triangulation
between radius $r$ and $r+1$ in the finite faces. The shaded areas are some
of the faces of the ball of radius $r+1$. The infinite face may contain
additional sub-faces.

This gives rise to a tree-like structure for the triangulation, as in
Figure~\ref{fig:ball}(b). Each outer face of the ball of radius $r$
corresponds to a vertex in the $r$'th level of the tree. The face
corresponding to a child is contained in the face of the parent vertex. An
infinite triangulation will yield an infinite tree. Similarly, if a
triangulation is one ended, then so is the corresponding tree, i.e., the
tree is composed of a single infinite branch from the root  with finite
sub-trees growing from it. Note that while any triangulation determines a
tree, the converse is false. The tree does not determine the triangulation.

A vertex in the tree corresponds to an external face of some triangulation,
so there are different types of vertices depending on the face
sizes. Labeling each vertex with the boundary size of the corresponding
face, we see that the UIPT gives rise to a multi-type tree process.

From Theorem~\ref{thm:locality} we see that
if we condition the first $r$ levels of the tree and on which vertex in the
$r$'th level is in the infinite branch, then Theorem~\ref{thm:locality}
tells us that the remaining sub-trees are independent. Thus, we see that at
each level, one vertex, with a known distribution, has an infinite sub-tree
above it, and the others have independent numbers of offspring of
independent types. The tree process is thus just a multi-type Galton Watson
process conditioned to survive. Without the conditioning we get the tree
corresponding to a free triangulation of the sphere, which we know to be
a.s.\ finite. However, since the free process has a power tail on its
size, it is critical. Thus, the above description is just the construction
of a critical Galton Watson process conditioned on survival (see
\cite{Lyons}).

\section{Type Relations}  \label{sec:equiv}     %\input equivalent.tex

The two types of UIPT are part of a wider class of random planar object
satisfying common properties. This was first hinted at by the universality
of the asymptotic formulas for counting various planar objects. Between
type II and type III triangulations there is a more fundamental relation,
enabling us to find a direct transformation between type II and type III
triangulations. A similar transformation also holds between type I and type
II triangulations as well as other pairs of classes of planar objects.

Roughly, the idea for passing from a type II triangulation to a type III
triangulation is to take each double edge and to remove all the triangles
and vertices inside it. The two edges are then glued together to get again
a triangulation of the plane or sphere as the case may be. Conversely, to
get from a type III triangulation to a type II one, we will take each edge
and replace it with a double (or multiple) edge with some distribution on
the triangulation inside the resulting 2-gons. Recall that we allowed the
triangulation of the 2-gon with no internal vertices, and gluing it in a
2-gonal outer face meant gluing the two edges together. Thus, with some
probability ($1/Z_0=8/9$ actually) this empty triangulation is used and the edge
remains a single edge.

Both directions pose some difficulties. A 2-gon partitions a triangulation
to two components. How do we decide which is the inside and which the
outside? In an infinite triangulation of the plane we wish to contract the
finite side, but for a finite triangulation of the sphere it is not so
clear. Also, there is the possibility that the root of the triangulation is
deleted in this way, and then a new root is needed.

In the opposite direction, there is the question of the distribution for
the triangulation of the 2-gon added. The natural candidate in the infinite
case is the free triangulation of the 2-gon. Again, in the finite setting
things are more delicate. Since we define the infinite triangulation
measures as limits of the finite ones, we need to find some transformation
of the finite measures first.

\medskip

For a type II triangulation $T$, when we contract 2-gons as described above,
until no double edges remain, the result is a maximal (with respect to
inclusion) 3-connected sub-triangulation, since the only way 2 vertices could
separate the graph is by forming a 2-gon. Therefore, the transition can be
summarized as taking a single maximal 3-connected subgraph of the
triangulation. The natural choice for this is to take the 3-connected
component containing the root triangle. This also saves us the trouble of
choosing a new root in the case that the old root was in one of the
contracted 2-gons. Note that the root vertex or even the root edge is not
enough, since 2 vertices may be in the intersection of two distinct
3-connected components. However the 3 vertices of any triangle determine a 
unique 3-connected component.

\begin{defn}
Let $T$ be a rooted type II triangulation.  Define $\tilde T$ to be the type III
triangulation composed of the 3-connected component of the root in $T$,
with all double edges identified into single edges. For a measure $\nu$ on
rooted type II triangulations let $\tilde \nu$ be the resulting measure on
type III triangulations, i.e., for any event $R$:
\[
\tilde \nu(R) = \nu\bigl (\{T\mid\tilde T \in R\}\bigr) .
\]
\end{defn}

This operation is known as taking the core of a structure \cite{BFSS}. In
general, for two classes of rooted planar objects, one more restricted than
the other, the core of a member of the wider class is its largest partial
structure containing the root included in the smaller class
(when it is unique).

\begin{lemma}\label{lem:transform}
For any finite $N$, for some coefficients $a_{n,i}$:
\[
\widetilde{\tau^2_N} = \sum_{i\leq N} a_{N,i} \tau^3_i .
\]
In the limit, for some constants $a_i,a_\infty >0$:
\[
\widetilde{\tau^2} = a_\infty \tau^3 + \sum a_i \tau^3_i .
\]
\end{lemma}

In the infinite case this means that the 3-connected component of the root
is either a finite sphere triangulation with some distribution on the size
where all triangulations of the same size are equiprobable, or it is an
infinite triangulation. Conditioned on the latter case it is just the
infinite type III UIPT. The asymptotics of the coefficients $a_{n,i}$ are
described in \cite{BFSS}

\begin{proof}
Consider first the finite case. All we need to show is that any two type
III triangulations have the same probability of appearing as the
3-connected component of the root in $\tau^2_N$, i.e., that for any
triangulation $U$ the number of triangulations $T$ with $|T|=N$ and
$\tilde T = U$ depends only on $|U|$.

This is clear, since any two triangulations of the same size have the same
number of edges. Specifically, $U$ has $3|U|-6$ edges. Formally, if
\mbox{$Z^2_0(x)= \sum_n \phi^2_{0,n} x^n$} is the generating function for
triangulations of a 2-gon, then the number of ways $U$ can come about is
the coefficient of $x^{N-|U|}$ in $(Z_0(x))^{3|U|-6}$, which is, of course,
determined by $|U|$.

\medskip

The infinite case follows from the finite case by taking a weak limit. 
The map $T \to \tilde T$ is continuous with respect to the topology on the spaces of
type II and III triangulations. Since $\tau^3_N$ is supported on
triangulations with $N$ vertices, they have disjoint supports for distinct
$N$. Therefore, necessarily:
\begin{eqnarray*}
\widetilde{\tau^2}
  &=& \lim_{N \to \infty} \widetilde{\tau^2_N}      \\
  &=& \lim_{N \to \infty} \sum_i a_{N,i} \tau^3_i       \\
  &=& \sum a_i \tau^3_i + a_\infty \tau^3 ,
\end{eqnarray*}
where $a_i = \lim_N a_{N,i}$ must exist and 
\[
a_\infty = \lim_{s \to \infty} \lim_{N \to \infty} \sum_{i>s} a_{N,i}
\]
is the part of the measure that tends to infinity. In the infinite case, we
can also give an explicit formula for $a_i$. This is done in much the same
way that we calculated the probability of a given ball when proving the
limit of $\tau_N$ exists. Indeed, to find the probability
$\widetilde{\tau^2}(T)$ for some type III triangulation $T$ with $|T|=n$ we
just need to find the probability $\tau^2(T)$ when each edge is replaced by
an external face of size 2. By Proposition~\ref{prop:sub_prob} this is:
\[
(3n-6) \alpha^{3-n} (Z^2_0)^{3n-7} C^2_0 / C^2_1 .
\]
(A sphere triangulation with $n$ vertices has $3n-6$ edges). Substituting
$Z^2_0=9/8$ and the values of $C$ this translates to:
\[
\frac{2^{19}}{3^7}(n-2)(\frac{27}{256})^n .
\]
Since there are $\phi^3_{n-3,1}$ possible triangulations of size $n$, the
probability $a_n = \tau^2(|\tilde T|=n)$ is:
\[
a_n=\frac{2^{20}(4n-11)!}{3^7(n-3)!(3n-7)!} (\frac{27}{256})^n .
\]
In order to find $a_\infty$ we need to sum $a_n$. Since $(256/27)^n\,a_n$
is a linear
combination of $\phi_{n-3,1}$ and $n\phi_{n-3,1}$, the generating function
$A(t) = \sum a_n (256 t/27)^n$ is a linear combination of $t^{-3}Z_1(t)$ 
and its derivative. Using that we find: $\sum a_n = 1/2$ and the remainder:
\[
a_\infty=1/2 .
\]
\end{proof}

Since we know that the UIPT is a.s.\ one ended, it has at most one infinite
3-connected component. The above calculations tell us more. We see that
with probability 1/2 the root triangle is part of the infinite 3-connected
component. In fact, if the root is in a finite 3-connected component, then
this component has a number of 2-gonal external faces, and the infinite one
contains a triangulation with the same LAW as the original UIPT. Iterating
this we see that there is always a unique infinite 3-connected component,
and with probability 1/2 the root is part of it. This is an infinite
version of an asymptotic result on finite triangulations found in
\cite{GaWo}. In fact, we know now the distribution of the size of the
3-connected component of the root, as $a_n$ is the asymptotic probability
that the component has size $n$.

How do we get back from the type III UIPT to the type II UIPT? We need to
find the distribution of a UIPT conditioned on including an infinite
triangulation. Theorem~\ref{thm:locality} deals with the UIPT conditioned
on containing a finite sub-triangulation, and by conditioning on a growing
subsequence of triangulations, we see that to get back from the infinite 
3-connected component to the whole type II triangulation we need to replace
each edge of $\tilde T$ with a free triangulation of a 2-gon.

Note that the expected number of triangles in a free triangulations of a
2-gon is twice the expected number of internal vertices and so is 2/3
(again, this is the derivative of $Z^3_0(t)$ at $\alpha^{-1}$). Since a 
triangulation contains 3/2 times as many edges as triangles, we see that in
some sense in the resulting type II triangulation 1/2 the triangles were in
the original type III triangulation and 1/2 were added.

\medskip

As a consequence of this relation, some results on the type II UIPT are
valid for type III as well. Those include the results on growth and on
percolation derived in \cite{UIPT2}, among others.

\bigskip

\filbreak
\begingroup
\small
\parindent=0pt

\vtop{  %
\hsize=2.3in  %
Omer Angel\\
Department of Mathematics\\  %
Weizmann Institute of science\\
Rehovot, 76100, Israel\\
{omer@math.weizmann.ac.il}
}
\bigskip  %
\vtop{  %
\hsize=2.3in  %
Oded Schramm\\  %
Microsoft Corporation\\  %
One Microsoft Way\\  %
Redmond, WA 98052, USA\\  %
{schramm@microsoft.com}  %
}  %
\endgroup  %

\filbreak

\end{document}